\documentclass[a4paper,12pt]{article}
\usepackage{amsthm,amsfonts,amsmath,amssymb}
\usepackage[utf8]{inputenc}
\usepackage[english]{babel}
\usepackage[final]{graphicx}
\usepackage{commath}
\usepackage{enumitem}
\usepackage{amsmath}

\title{Disjointly almost trivial unbounded functionals}
\author{Konstantin Storozhuk}
\date{December 2025}

\begin{document}

\maketitle

Let $X$ be a Banach lattice $l_p(\Bbb N)$.
We show that there exist unbounded functionals on $X$ that take at most one nonzero value on an arbitrary family of elements whose supports are pairwise disjoint. This result, in particular, provides a negative answer to questions 1–6 of the paper \cite{1}, which ask whether certain "goodness" properties of a linear operator on disjoint families imply "goodness" properties of the operator as a whole.

\vskip1mm

The construction is based on the following fact:\vskip1mm

\bf Theorem 1\it. Let X =$l_p(\Bbb N)$. There exists a vector subspace $V$ of $X$ of infinite codimension such that for an arbitrary disjoint family $A\subset X$ all elements of $A$, with the possible exception of one, belong to $V$.\vskip1mm

Proof\rm. 
Let $\mathcal{F}$ be a free ultrafilter on $\Bbb N$.
We can take $V$ to be the set of all functions $x$ such that $0(x)\in \mathcal{F}$. Here $0(x)\subset \Bbb N$ is the set on which $x$ is zero.
Equivalently, $x\in V$ if and only if $supp(x)\notin \mathcal{F}$.

 $V_\mathcal{F}$ is a vector subspace. This follows from the fact that the filter is closed under intersection and superset operations. For example, if $x,y\in V_\mathcal{F}$, then $x+y\in V_\mathcal{F}$, since $0(x)\cap 0(y)\subset 0(x+y)$. Further, if $a$ and $b$ are disjoint, then $0(a)\cup 0(b) =\Bbb N$, which implies that at least one of these sets belongs to $\mathcal{F}$ (here we use the ultrafilter property: for an arbitrary subset $A \subset
M$, exactly one of the sets $A$, $X\backslash A$, belongs to $\mathcal{F}$). That the space $V_\mathcal{F}$ has infinite codimension is easy to deduce, for example, from the fact that $V_\mathcal{F}$ does not contain sequences $x$ for which the set $0(x)$ is finite. Let us show that the codimension of the constructed space $V$ is not less than the cardinality of the continuum. The codimension of the space $V_\mathcal{F}$ is uncountable. Indeed, for each number $\alpha<0$ we can consider the sequence $x^\alpha:=(2^{n\alpha})$. It is easy to see that these sequences are linearly independent; moreover, no finite nontrivial linear combination of them lies in the space $V_\mathcal{F}$. Indeed, any finite such combination $w$ of sequences with exponents $\alpha_1<\alpha_2<\cdots<\alpha_{max}<0$ has an asymptotics of the form $cx^{\alpha_{max}}$ and, therefore, starting from some number, is always nonzero. Therefore, the set $0(w)$ is finite and $ 0(w)\notin\mathcal{F}$. The theorem is proved.\vskip1mm

\it Remark 1\rm. The subspace $V_\mathcal{F}$ contains all elements whose support is finite, and therefore this space is dense in $X$ if $p<\infty$.

\it Remark 2\rm. The arguments in the theorem remain valid for function spaces on an arbitrary infinite set.

\it Remark 3\rm. If the ultrafilter is not free, then it will consist of supersets of some singleton $\{n_0\}$. In this case, our construction will lead to the construction of a space consisting of sequences whose coordinate with $n_0$ is zero. This space is closed, and its codimension is $1$.

\vskip1mm

\bf Theorem 2\it. Let $X=l_p$, $1\leq p\leq \infty$. There exists a unbounded linear functional $h:X\to \Bbb R$ such that for any disjoint family of elements of $X$, all these elements, except possibly one, lie in $\ker h$.\vskip1mm

Proof\rm. Let $V_\mathcal{F}$ be a subspace of $X$ constructed as in Theorem 1. The linear functional extension theorem and the infinite codimension of $V_\mathcal{F}$ allow us to extend the zero functional defined on $V_\mathcal{F}$ to a unbounded functional on all of $X$. The theorem is proved.

Note that the construction of the functional $h$ is "doubly non-constructive": we use both the existence of free ultrafilters and the linear  extension theorem.

\vskip1mm

\it Remark 4\rm. If $p<\infty$, then \it{any} \rm nontrivial linear functional on $l_p$ whose kernel contains $V_\mathcal{F}$ is discontinuous due to the density of $V_\mathcal{F}$ in $X$. At the same time, on the space $l_\infty$ there already exist nontrivial bounded functionals that are equal to zero on $V_\mathcal{F}$. For example, such is the functional $(\mathcal{F}$-$\lim)$, which associates to each bounded sequence its limit over an ultrafilter $\mathcal{F}$ .


\vskip3mm




\begin{thebibliography}{999}

\bibitem{1}
E. Yu. Emelyanov, N. Erkurşun-Özcan, S. G. Gorokhova,  $d$-Operators in Banach Lattices.
Siberian Mathematical Journal,
V. 66, p. 1499–1508 (2025).



\end{thebibliography}
\end{document}